\documentclass{article}
\usepackage[utf8]{inputenc}
\usepackage{tikz}
\usetikzlibrary{graphs} 
\usepackage{amsmath, amssymb, amsthm}
\usepackage{graphicx}
\usepackage{hyperref}
\usepackage{iftex}

\ifpdf
  \usepackage{underscore}         
  \usepackage[T1]{fontenc}        
\else
  \usepackage{breakurl}           
\fi

  \newtheorem{theorem}{Theorem}[section]
  \newtheorem{lemma}{Lemma}[section]
  
  \newtheorem*{theorem*}{Theorem}

  \newtheorem{definition}{Definition}[section]

\title{Proof verification by polynomial fingerprinting\footnote{A preliminary version was presented at the Working Formal Methods Symposium FROM 2025, Iași, România, and was published in the proceedings Electronic Proceedings of Theoretical Computer Science EPTCS 427 as {\it Polynomial Fingerprtinting for Trees and Formulas}, see also arXiv:2509.11877. The present version contains better and more consecvent notations, the true origin of the problem,  some more explanatins and more details for first-order logic, respectively for other logic systems.}}
\author{Mihai Prunescu
\thanks{(1) Research Center for Logic, Optimization and Security (LOS), Faculty of Mathematics and Computer Science,  University of Bucharest, Academiei 14, 010014 Bucharest, Romania. mihai.prunescu@gmail.com}
\thanks{(2) Simion Stoilow Institute of Mathematics of the Romanian Academy,  Research unit 5, P. O. Box 1-764, RO-014700 Bucharest, Romania. mihai.prunescu@imar.ro}
}

\begin{document}

\maketitle

\begin{abstract}
To cater to the needs of fast verification for mathematical proofs, we describe a method to encode formal sentences in $2 \times 2$ - matrices over multivariate polynomials with integer coefficients. This correspondence is homomorphic:  usual proof-steps like modus-ponens or variable substitution in terms and formulae become operations with matrices. By evaluating the polynomial  variables in random elements of a suitably chosen finite field, the proof is replaced by a numeric sequence. Only the values corresponding to axioms and tautologies have to be computed from scratch. The values corresponding to derived formulas are computed from the values corresponding to their ancestors by applying the homomorphic properties. The polynomial matrix corresponding to the conclusion of the proof is also evaluated in the chosen random values. If the last term of the numeric sequence equals the evaluation of the conclusion, by the Schwartz-Zippel Lemma, the proof is with high probability correct. 


\end{abstract}

\section{Motivation}  \label{pf:motivation} 

 There is an increasing interest in authomatic proof checkers, like  Methamath,  Hol Light Coq,  Isabelle, etc. Moreover, there are entire research programs to rewrite mathematics on automatically verifiable bases. On the other hand, there is an increasing economical and social interest in correctness certificates for programs or for program runs, like online banking or cryptocurrency transactions. The computation speed is an important factor for both directions of development, but is more sensible in economic applications. The user needing a correctness certificate for a transaction with cryptocurrencies would accept also a probabilistic verification with big probability, if this verification is fast enough to be integrated in his business flow. Differently, the fundamental mathematician would prefer a verification with $100\%$ surety, even if it takes a little bit longer. Remarkably, a correctness certificate for the result of a program can be also given in the form of a mathematical proof, see the K semantic framework \cite{RosuSerbanuta}. The correctness proofs can be done in matching logic, see \cite{matching1} and \cite{matching2}.  The first idea how to secure a mathematical proof by some zero-knowledge procedures goes back to Blum \cite{Blum}.

There is progress in using Zero Knowledge methods for generating correctness certificates for some special formats of mathematical proofs, see Couillard et al. \cite{Couillard} or Luick et al. \cite{Luick}. Dominant trends in runtime verification are bazed on Zero Knowledge Virtual Machines based on SNARK or STARK algorithms, like for example Risk0 \cite{Risc0}, but all the procedures implemented so far are very slow, as they encode all the hardware procedures which took place while the program was running. From this point of view, the approach over a verifiable mathematical proof generated over the K semantic framework looks to be more robust. Unfortunately, the $100\%$ verification of a mathematical proof is again very slow. 

The goal of this paper is to sketch a fast verification procedure for mathematical proofs, which certifies its correctness with high probability. 
To this end, one finds an encryption of formalized statements in polynomial matrices, such that proof steps like modus-ponens and substitution are homomorphically translated in algebraic operations with matrices. This will allow a fast verification procedure which works according to the Schwartz-Zippel Lemma \cite{Schwartz}, \cite{Zippel}. 

As the proof verification is finally done by evaluating an arithmetic circuit,  it can be combined with Zero Knowledge procedure for arithmetic circuits, like CIRCOM \cite{Circom} and with the Fiat-Shamir method to transform interactive protocols in correctness certificates \cite{FiatShamir}, to be used to produce Zero Knowledge correctness certificates for mathematical proofs.

\section{The problem}\label{pf:problem}

There are already some methods of fingerprinting, designed for proof verification, and finally for the generation of succint certificates via CIRCOM, as communicated to the author by Brandon Moore, \cite{Moore}. In all these methods, proof-lines are considered to be strings. With strings, it is quite easy to associate numbers, respectively elements of large finite fields, such that the modus ponens step works homomorphically. But it is more difficult to find a homomorphic correspondence between strings and field-elements with respect to the substitution of a variable with a whole formula or with a term. More exactly, the homomorphic substitution of strings, which works more or less like the positional notation for numbers, is homomorphic only for one occurrence of a variable with a term, because for any further occurrence the position of the variables change after the substitution of the first occurrence. Suppose that the free variable $x$ occurs in $\varphi$ say $15$ times, and should be replaced with the term $t$ to build the result of this substitution, usualy denoted by $\varphi(x/t)$. Then for every occurrence, one needs a proof-line. Even worse, in some approaches one needs supplementary proof-lines to compute the positions of the remaining occurrences.

The approach presented in this article starts with formulas considered as trees, instead of strings. It has the disadvantage that for every function symbol or relation symbol of arity $k$ it introduces $k+1$ variables. Also, a faithful translation of formulas to polynomials, respectively to elements of large finite fields, needs the use of some essentially non-commutative algebras. However, this method allows both modus-ponens and substitution to be done homomorphically. The substitution can be done in one step, independently of the number of occurrences of the variable in the formula. Also, the method allows extensions like use of quantifiers, use of modal operators, application of patterns or some other operators, like the least fix-point operator. 

The present approach is not meant to replace other approaches. It is only meant to complete them, by making the substitution steps faster, while the modus ponens steps remain homomorphic. 

An outline of this paper: In Section \ref{pf:motivation} is presented some motivation, based on general needs of the human society. In Section \ref{pf:problem}, the motivation became technical. Section \ref{pf:principle} shows the solution in principle, as it looks from the outside. In this section is also shown that the Schwarz-Zippel Lemma assures that positive verification over polynomial evaluation means correctness of the proof with big probability. In Section \ref{storysimp} we introduce the non-commutative ring of matrices $M_{2 \times 2}(\mathbb Z[x_1, x_2, \dots, x_n, \dots])$. For a polynomial variable $x$, a matrix $A(x)$ is introduced. This matrices have the property that a product $A(y_1)A(y_2) \dots A(y_m)$, when $y_1, \dots, y_m$ belongs to the set of variables, possibly with repetitions, uniquely determines the number of variables and their identity. This algebra is used in Section \ref{treesandpolynomials} to construct the unique polynomial matrix associated with a formula. In Section \ref{homomorphprop} is shown why modus ponens and substitution are homomorphic operations with respect with the association announce above. The proof is done for an unspecified logic. Section \ref{example} contains a worked out example of formalized proof in propositional logic, which has been fingerprinted. Section \ref{firstorder} presents more details about fingerprinting formulas and proofs of first order logic, and some tricks and recommendations for formulas and proofs in some further proof-systems.

\section{The principle}\label{pf:principle}

Suppose that one fixed a formalism for mathematical proofs, which includes propositional and first order logic, but may include also other particular proof systems, like matching logic \cite{matching1}, \cite{matching2}. One considers  a subsystem consisting of a finite signature and a finite number of axioms and axiom schemes.  To every well formed expression $E$, which can be a term or a formula, one associates a matrix $[E] \in M_{2 \times 2}(\mathbb Z[x_1, x_2, \dots, x_n, \dots])$, such that a matrix corresponds to at most one formula. 

The proofs in our system are structured as follows. Every proof-line starts with an explanation about the origin of the statement. This may be: axiom $5$, axiom $3$ for variables $a$, $b$, $c$, modus ponens of line $7$ with line $5$, or line $8$ with all occurrences of variable $d$ replace by term (or formula) $t$, which has to be given explicitly. The explanation is followed by the resulting statement. 

This association between expressions and matrices has homomorphic properties in the following sense. For every deduction rule, there is a way to compute the matrix corresponding to the conclusion from the matrices corresponding to the premises. For example, suppose that the formula $\psi$ is obtained by applying modus ponens to premises $\varphi$ and $\varphi \rightarrow \psi$. Then there is an algebraic term $MP(A, B)$ such that always:
$$ [\psi] = MP([\varphi], [\varphi \rightarrow \psi]).$$
Also, if an expression $F$ results from the substitution of a term (formula) $t$ for all occurrences of the variable $a$ in the expression $E$, there is an algebraic term $Subst_a(A, B)$ such that:
$$[F] = Subst_a([E], [t]).$$
It will come out that this step needs some auxiliary matrices to be already computed, because the matrix corresponding to the substitution has the form:
$$[F] = [E] - [E]_a \cdot [a] + [E]_a \cdot [t]. $$ 
Here, $[a]$ is the polynomial matrix corresponding to the proof variable $a$, which might mean a variable of propositional logic, an element variable in first order logic, or an element variable of a given sort in matching logic. One should carefully make the difference between variables in the sense of the mathematical proof and polynomial variables used to encode the proof in polynomial matrices. The matrix $[E]_a$ depends on the variable $a$ and on the well formed expression $E$ and should be already computed in the moment of substitution. 

Given this necessity, one associates with an expression $E$ containing free proof variables (logical variables or first-order element variables) $a_1$, $\dots$, $a_k$ a tuple of matrices:
$$P(E) = ([E], [E]_{a_1}, \dots, [E]_{a_k}),$$ 
the {\bf algebraic fingerprint} of $E$. 

Once this machinery works, one can use it for fast proof verification as follows:
\begin{enumerate}
\item Choose a large finite field $\mathbb F$.
\item Suppose that the polynomial variables $x_1, \dots, x_n$ arise in the polynomial matrices $F(E)$, where $E$ are various well formed expressions present in the mathematical proof. One chooses random elements $r_1, \dots, r_n \in \mathbb F$. 
\item One constructs all matrices corresponding to proof lines (proof statements) using elements $r_i$ instead of polynomial variables $x_i$. This is equivalent to evaluate the corresponding polynomial matrix in this values. For any expression $E$, one denotes by
$$[[E]] = [E](r_1, \dots, r_n)$$
the {\bf numeric fingerprint} of $E$.
\item The other components of the algebraic fingerprint are also evaluated as
$$[[E]]_{a_s} = [E]_{a_s}(r_1, \dots, r_n).$$
\item One computes the numeric fingerprint of the conclusion $\gamma$ inductively, following all proof-lines. For the lines representing instances of the axioms, the corresponding fingerprints are computed by direct translation. For the lines obtained by modus ponens or by substitution, the corresponding fingerprints are computed using the homomorphic properties. We obtain the fingerprint
$$[[\gamma]]'$$
as result of the last deduction, which could be a modus ponens, a substitution, writing or eliminating a quantifier, or another step which can be done by a homomorphic rule. 
\item Given the formal statement $\gamma$, which has to be the conclusion of the mathematical proof, one directly computes its numeric fingerprint by definition (direct translation):
$$[[\gamma]]'' = [\gamma](r_1, \dots, r_n).$$
\item One compares the matrices $[[\gamma]]'$ and $[[\gamma]]''$. If they are different, the proof is not correct. If they are equal, the proof is correct with high probability.
\end{enumerate}

This verification procedure works grace to the Theorem of Schwartz \cite{Schwartz} and Zippel \cite{Zippel}:

\begin{theorem}
    Let $\mathbb F$ be a finite field and let $f \in \mathbb F[x_1, \dots, x_n]$ be a non-zero polynomial of degree $d \geq 0$. If $r_1, r_2, \dots, r_n$ are selected randomly and if the choices are independent in $\mathbb F$, then:
    $$ Pr[f(r_1, r_2, \dots, r_n) = 0] \leq \frac{d}{|\mathbb F|}.$$
\end{theorem} 

The Schwartz and Zippel Theorem  says that for a polynomial which is not identically zero, the probability that a random evaluation is zero in a big finite field is reasonably small. Let again $\gamma$ be the conclusion of the mathematical proof. Let $[\gamma]'' \in M_{2 \times 2} (\mathbb Z [x_1,x_2, \dots]) $ be the polynomial matrix obtained from $\gamma$ by direct encoding. Let $[\gamma]' \in  M_{2 \times 2} (\mathbb Z [x_1, x_2, \dots])$ be the polynomial matrix obtained from the encoding of the axioms and using the proof steps and their homomorphic properties. One proves with high confidence that $[\gamma]' = [\gamma]''$. To this scope, one computes their evaluations $[[\gamma]]'$ and $[[\gamma]]''$, where $[[\gamma]]'$ is computed from evaluations of the axioms and using homomorphic properties, and $[[\gamma]]''$ is computed by direct translation, but using only numeric matrices, instead of polynomial matrices. If the field is sufficiently large, the equality $[[\gamma]]'=[[\gamma]]''$ means that, with high probability, $[\gamma]' = [\gamma]''$, so the formula resulted from the proof is indeed identical with the claimed conclusion, which has been encoded directly. One can repeat the procedure with other random field elements in order to increase the requested level of confidence. 

More exactly, the matrix $[\gamma]'$ consists of following polynomials:
$$ [\gamma]' = \begin{pmatrix}
Q'(x_1, \dots, x_x) & P'(x_1, \dots, x_n) \\ 0 & k'
\end{pmatrix},$$
while the matrix $[\gamma]''$ looks like:
$$ [\gamma]'' = \begin{pmatrix}
Q''(x_1, \dots, x_x) & P''(x_1, \dots, x_n) \\ 0 & k''
\end{pmatrix},$$
where $k'$ and $k''$ are integers, and under the evaluation in a field $K = \mathbb F_p$, where $p$ is a large prime, go to $k' \bmod p$, respectively $k'' \bmod p$.
The proof is correct if and only if the proven formula, uniquely represented by $ [\gamma]'$, is identic with the formula to prove, uniquely represented by $ [\gamma]''$. But $ [\gamma]'$ and $ [\gamma]''$ are not computed. Via polynomial evaluation in $(x_1, \dots, x_n) = (r_1, \dots, r_n) \in K^n$, the equality $[[\gamma]]' = [[\gamma]]''$ implies that $(Q' - Q'')(r_1, \dots, r_n) = 0$, $(P' - P'')(r_1, \dots, r_n) = 0$ and 
$(k'-k'') \bmod p = 0$. This clearly means that:
\begin{enumerate}
\item If $[[\gamma]]' \neq [[\gamma]]''$, then the proof is false with probability $1$.
\item If $[[\gamma]]' = [[\gamma]]''$, then the proof is correct with probability $\leq n/p$, where $n$ is the number of polynomial variables occurring in the fingerprints.
\end{enumerate}
We will see that $n$ is strictly less than the number of constants $+$ the number of element- and propositional variables $+$ the number of non-logic symbols times the maximal arity of the non-logic symbols, increased by one. In fact, for every function symbol or relation symbol of arity $k$, one adds $k+1$ polynomial variables. In conclusion, the number of polynomial variables is relatively small and completely previsible, so there are no difficulties to choose a prime $p$ which is large enough to assure a large correctness probability in just one check.

\section{Matrices of multivariate polynomials}\label{storysimp} 

The original observation which led to this subsection was that matrices consisting of different variables:
$$ M(k) = \begin{pmatrix}
    x_{4k+1} & x_{4k+2} \\ x_{4k+3} & x_{4k+4}
\end{pmatrix}$$
are non-commutative to such extent that if two products are equal: 
$$M(i_1)\dots M(i_n) = M(j_1) \dots M(j_m)$$
then $n = m$ and $i_1 = j_1$, $\dots$, $i_n = j_n$. This means that such a monomial (product of elementary matrices) contains information about the number of factors, their order, and their identities. All three elements of information are essential to encode a path inside a tree.

However, proceeding with such matrices would be expensive because one has to choose four field elements to evaluate every elementary matrix and to keep four elements for every matrix to be kept as part of a fingerprint. Instead, one presents a system of elementary matrices that achieves the same goal, but needs just one field element to be evaluated and only three field elements for every matrix which is part of a fingerprint. 

\begin{definition}\label{defmatrixvar}
    Let $x_1$ be a polynomial variable. Let:
    $$A(x_1) = \begin{pmatrix}
        x_1 & 1 \\ 0 & 1
    \end{pmatrix}$$
    Observe  that $A(x_1) \in M_{2 \times 2}(\mathbb Z[x_1, x_2, \dots])$. 
\end{definition} 

We show now that these elementary matrices fulfill the same property: a monomial contains enough information to uniquely determine the number of elementary matrices, their order and their identities.

\begin{lemma}\label{lemmamatrix}
    Consider a set of different variables $V = \{x_1, x_2, \dots, x_k\}$. Suppose that $0 \leq i_1, \dots, i_n, j_1, \dots, j_m \leq k$. If:
    $$A(x_{i_1}) A(x_{i_2}) \dots A(x_{i_n}) = A(x_{j_1}) A(x_{j_2}) \dots A(x_{j_m}),$$
    then the following equalities take place: $n = m$, $i_1 = j_1$, $\dots$, $i_n = j_n$.
\end{lemma}

\begin{proof}
    If in this identity, we set $x_1 = \dots = x_k = 1$, as we observe that:
    $$\begin{pmatrix}
        1 & n-1 \\ 0 & 1
    \end{pmatrix}
    \begin{pmatrix}
         1 & 1 \\ 0 & 1
    \end{pmatrix} = 
    \begin{pmatrix}
        1 & n \\ 0 & 1
    \end{pmatrix},$$
    we get that $n = m$. 
    Let us denote with $S(n)$ the statement: {\it If
     $$A(x_{i_1}) A(x_{i_2}) \dots A(x_{i_n}) = A(x_{j_1}) A(x_{j_2}) \dots A(x_{j_n})$$
    then $i_1 = j_1$, $\dots$, $i_n = j_n$.}
    We observe that $S(1)$ is evident by identifying the entries. Suppose that we have proved $S(n)$. We look at the hypothesis of $S(n+1)$:
    $$A(x_{i_1}) A(x_{i_2}) \dots A(x_{i_{n+1}}) = A(x_{j_1}) A(x_{j_2}) \dots A(x_{j_{n+1}}).$$
    Observe by induction that:
    $$A(x_{i_1}) A(x_{i_2}) \dots A(x_{i_n}) = 
    \begin{pmatrix}
        x_{i_1}x_{i_2} \dots x_{i_n} & P(x_{i_1}, x_{i_2}, \dots, x_{i_{n-1}}) \\ 0 & 1
    \end{pmatrix},
    $$ where $P(z_1, \dots, z_{n-1}) \in \mathbb Z[z_1, \dots, z_{n-1}]$ is a fixed polynomial, {\it with the property that no variable $z_i$ divides $P$}. We write the hypothesis of induction in the form:
    $$
    \begin{pmatrix}
        x_{i_1} & 1 \\ 0 & 1
    \end{pmatrix}
    \begin{pmatrix}
        x_{i_2}x_{i_3} \dots x_{i_{n+1}} & P(x_{i_2}, x_{i_3}, \dots, x_{i_{n}}) \\ 0 & 1
    \end{pmatrix} = $$$$ =
    \begin{pmatrix}
        x_{j_1} & 1 \\ 0 & 1
    \end{pmatrix}
    \begin{pmatrix}
        x_{j_2}x_{j_3} \dots x_{j_{n+1}} & P(x_{j_2}, x_{j_3}, \dots, x_{j_{n}}) \\ 0 & 1
    \end{pmatrix},
    $$ 
     so:
     $$
      \begin{pmatrix}
        x_{i_1} x_{i_2}x_{i_3} \dots x_{i_{n+1}} & 1 + x_{i_1} P(x_{i_2}, x_{i_3}, \dots, x_{i_{n}}) \\ 0 & 1
    \end{pmatrix} = $$$$ =
      \begin{pmatrix}
        x_{j_1} x_{j_2}x_{j_3} \dots x_{j_{n+1}} & 1 + x_{j_1} P(x_{j_2}, x_{j_3}, \dots, x_{j_{n}}) \\ 0 & 1
    \end{pmatrix}.
     $$
     Identify the corresponding entries:
     \begin{eqnarray*}
          x_{i_1} x_{i_2}x_{i_3} \dots x_{i_{n+1}} &=& 
           x_{j_1} x_{j_2}x_{j_3} \dots x_{j_{n+1}}, \\
           1 + x_{i_1} P(x_{i_2}, x_{i_3}, \dots, x_{i_{n}}) &=&
           1 + x_{j_1} P(x_{j_2}, x_{j_3}, \dots, x_{j_{n}}).
     \end{eqnarray*} 
     Apply the property that no variable $z_i$ divides $P$. By variable identification we get $x_{i_1} = x_{j_1}$. Multiply the hypothesis with $A(x_{i_1})^{-1}$ from the left-hand side.  One gets an instance of $S(n)$ and we apply the induction hypothesis. Of course  $A(x_{i_1})^{-1}$ does not belong to the matrices over polynomials, but to the matrices over rational functions. It is important that one can simplify with $A(x_i)$. 
\end{proof}

\section{Non-commutative representation of edges and nodes}\label{treesandpolynomials} 

 The goal of this section is to show how to further associate a matrix to a term or a formula represented by trees using the previous construction.

 Below, by {\bf edge variable} we understand a polynomial variable $x_i$. To every edge of the tree, and to every vertex, we will associate  a matrix of the shape:

 $$ A(x_i) =  \begin{pmatrix}
     x_i & 1 \\ 0 & 1
 \end{pmatrix}.$$
 
In order to represent formulas by trees, both logical and term-building operations are represented as vertices of the tree. For every specific symbol $c$ of arity $d=d(c)$ a number of $d + 1$ different fixed  edge variables $C, C_1, \dots, C_d \in \{x_1, x_2, \dots\}$ are associated. 

Suppose that a tree $T$ has root $c$ and the sub-trees connected with $c$ are $T_1, \dots, T_d$. Suppose that one already associated matrices $$[T_1], \dots, [T_d] \in M_{2 \times 2}(\mathbb Z[x_1, x_2, \dots])$$ with these sub-trees. Then we associate with $T$ the pair:
$$[T] = A(C) + A(C_1) [T_1] + \dots + A(C_d) [T_d],$$
where $C, C_1, \dots, C_d$ are the associated edge variables. 

\begin{definition}
    The height of a tree consisting of just one node is $0$. Consider the tree represented by the polynomial matrix: 
$$[T] = A(C) + A(C_1) [T_1] + \dots + A(C_d) [T_d].$$
If the trees $T_1$, $\dots$, $T_d$ have heights $t_1$, $\dots$, $t_n$ respectively, then the height of $T$ is 
$$1 + \max(t_1, t_2, \dots, t_d).$$
\end{definition}

\begin{definition}
    If $\varphi$ is a formula or a term, let $[\varphi]$ denote the polynomial matrix associated with its tree. 
\end{definition}

\begin{theorem}
    A matrix $A \in M_{2 \times 2}(\mathbb Z[x_1, x_2, \dots, x_n, \dots])$ represents at most one well formed expression (term or formula).
\end{theorem}

\begin{proof} We fix a logic signature.  The proof works by induction on the height of the tree representing the expression, called in the sequent {\bf height of the expression}.

 If the tree corresponding with an expression consists of only one vertex, its matrix is nothing but the matrix corresponding to the associated vertex variable and has the shape
$$A(x) = \begin{pmatrix} x & 1 \\ 0 & 1 \end{pmatrix}.$$
So we have proved the statement for trees of height $0$. Suppose that the statement has been already proved for all expressions of height at most $n$. We consider an expression of height $n+1$ represented by the matrix:
$$[T] = A(C) + A(C_1) [T_1] + \dots + A(C_d) [T_d].$$
Then the expressions corresponding to $T_1$, $\dots$, $T_d$ have heights $\leq n$, so by the induction hypothesis they are uniquely determined by these polynomial matrices, and these expressions have uniquely determined positions in the argument vector of the operation or relation symbol encoded by the vertex variable $A(C)$. These positions are given by the used edge-variables $A(C_1)$, $\dots$, $A(C_d)$. 
\end{proof}

{\bf Example}:  Consider the following inductive definition:
\begin{enumerate}
    \item The letters $x$, $y$, $z$ are atomic propositional formulas.
    \item If $\varphi$ and $\psi$ are formulas, then:
    $$\neg\,\varphi,\,\, \varphi \rightarrow \psi, $$
    are formulas.
\end{enumerate}
 
The alphabet is $A = \{x, y, z, \neg, \rightarrow \}$. 

The variables $x, y, z$ are symbols of arity $0$ and will always be final nodes. We associate them with the matrices:
$$
[x] =  A(X) = \begin{pmatrix}
    X & 1 \\ 0 & 1
\end{pmatrix},\,\,\,\,
[y] = A(Y) = \begin{pmatrix}
    Y & 1 \\ 0 & 1
\end{pmatrix},\,\,\,\,
[z] = A(Z) = \begin{pmatrix}
    Z & 1 \\ 0 & 1
\end{pmatrix}.
$$

The symbols with positive arity are $\{\neg, \rightarrow\}$. We associate with $\neg$ the matrices:
$$
A(N) = \begin{pmatrix}
    N & 1 \\ 0 & 1
\end{pmatrix},\,\,\,\,
A(N_1) = \begin{pmatrix}
    N_1 & 1 \\ 0 & 1
\end{pmatrix}.
$$
We associate with $\rightarrow$ the matrices:
$$
A(I) = \begin{pmatrix}
    I & 1 \\ 0 & 1
\end{pmatrix},\,\,\,\,
A(I_1) = \begin{pmatrix}
    I_1 & 1 \\ 0 & 1
\end{pmatrix},\,\,\,\,
A(I_2) = \begin{pmatrix}
    I_2 & 1 \\ 0 & 1
\end{pmatrix}.
$$
The $7$ variables $X, Y, Z, N, N_1, I, I_1, I_2$ are pairwise different.

The inductive steps are given by:
    $$[\neg \,\alpha] = A( N) + A(N_1) [\alpha],$$
    $$[\alpha \rightarrow \beta] = A(I) + A(I_1) [\alpha] + A(I_2) [\beta], $$
The statement of the Theorem is proved by induction over the building rules for formulas. What we really prove is the equivalent statement: {\it Every formula is encoded by only one matrix of polynomials}.  If $\varphi$ is an atomic propositional symbol, then $[\varphi]$ is $[x]$,  $[y]$ or $[z]$ and so from $[\varphi] = [\varphi']$ follows immediately $\varphi = \varphi'$. Suppose that $\varphi = \neg \, \alpha$. Then:
$$[\varphi] = A(N) + A(N_1) [\alpha].$$
We observe that $A(N)$ is the only one monomial of degree one present here. So one can conclude that we are reading a negation. All other monomials start with $A(N_1)$ because, as shown in Lemma \ref{lemmamatrix}, all these non-commutative monomials can start only with $A(N_1)$. Now, by the induction hypothesis, the formula $\alpha$ is uniquely encoded by $[\alpha]$, and it follows that $\varphi$ is uniquely encoded by $[\varphi]$. 

Now we consider the case $\varphi = \alpha \rightarrow \beta$. We have seen that: 
$$[\varphi] = A(I) + A(I_1) [\alpha] + A(I_2) [\beta].$$
Again $A(I)$ is the only one monomial of degree one and its presence shows that we are reading an implication. All other monomials have the shape $A(I_1) B$ or $A(I_2)B$. By the unicity of products of elementary matrices (Lemma \ref{lemmamatrix}), this monomials can start only with $A(I_1)$ or with $A(I_2)$. By common factor, we get the expression $A(I_1) [\alpha] + A(I_2) [\beta]$. As by induction hypothesis the formulas $\alpha$ and $\beta$ are uniquely expressed by the polynomial matrices $[\alpha]$, respectively $[\beta]$, it follows that $\varphi$ is uniquely expressed by the matrix of polynomials $[\varphi]$. 

\section{Homomorphic properties}\label{homomorphprop}

In this subsection we define the notion of fingerprint of a formula and we show that this notion enjoys homomorphic properties of the operations with formulas used in proofs: modus ponens and substitution. We show how the fingerprint of a formula produced by modus ponens can be computed from the fingerprints of the arguments of the operation modus ponens. Also we show how the fingerprint of a formula obtained by substitution can be computed from the fingerprints of the formula in which the variable has been substituted and the fingerprint of the formula (or term) which has substituted this variable. 

Suppose that three different lines of a proof read:
\begin{eqnarray*}
    \varphi \,\,\,\,\,\,\,\,\,\,\,\,\,\,\\
    \varphi \rightarrow \psi \\
    ---\\
    \psi
\end{eqnarray*}
such that the formula $\psi$ is deduced from the formulae $\varphi$ and $\varphi \rightarrow \psi$ by modus ponens. Suppose that we equip the implication symbol $\rightarrow$ with three matrices $A(I)$, $A(I_1)$ and $A(I_2)$ such that:
$$[\varphi \rightarrow \psi] = A(I) + A(I_1) [\varphi] + A(I_2)[\psi].$$
Then one can compute the conclusion as follows:
$$[\psi] = {A(I_2)}^{-1} \left ([\varphi \rightarrow \psi] - A(I) - A(I_1)[\varphi]\right ).  $$ 

Substitution also enjoys a homomorphic property. Suppose that one has a formula $\varphi(x)$ and substitutes $x$ with a tree $[\psi]$ corresponding to a formula or a term. We observe that:
$$[\varphi(x)] = \sum_{\text{nodes }c} A(X_{i_1})\dots A(X_{i_n}) \cdot A(X_c).$$
Here for every node $c$, the monomial $A(X_{i_1})\dots A(X_{i_n})$ consists of the edge-variables on the path from the root to the node $c$. If two such nodes are marked with $x$ and are to be substituted, one has:
$$[\varphi(x/\psi)] = [\varphi(x)] - A(X_{i_1})\dots A(X_{i_n})  [x] - A(X_{j_1})\dots A(X_{j_m}) [x] +$$ $$+ A(X_{i_1})\dots A(X_{i_n}) [\psi] + A(X_{j_1})\dots A(X_{j_m})[\psi].$$
In general, let $x$ be a propositional or a first-order element variable, and let $X$ be the polynomial variable associated to this symbol of arity $0$. Let $\varphi$ be a formula or a term. We denote by:
$$\sum _{c = x} A(X_{i_1})\dots A(X_{i_n}) \cdot A(X_c) := [\varphi]_x \cdot A(X_c).$$
It follows that in general for every formula or term $\psi$,
$$[\varphi(x/\psi)] = [\varphi] - [\varphi]_x \cdot A(X) + [\varphi]_x \cdot [\psi]. $$ 
Observe that the matrix $[\varphi]_x$ has been implicitly defined in the precedent formula. 

\begin{definition}\label{deffingerprint}
Let $\varphi$ be a well-formed expression over $A$, i.e. a term or a formula. Suppose that $x_1, \dots, x_k$ are the free variables in $\varphi$, which may be propositional variables or first-order element variables. We call the fingerprint of $\varphi$ the tuple:
$$([\varphi], [\varphi]_{x_1}, \dots, [\varphi]_{x_k}).$$
We denote the algebraic fingerprint of $\varphi$ with $F(\varphi)$.
\end{definition} 

Observe that, as we announced in Section \ref{pf:principle}, we will deal with two kinds of fingerprints. To a string $\varphi$ which might be a formula or a term, we have defined the fingerprint consisting of a vector of matrices over the polynomial ring $\mathbb Z [x_1, x_2, \dots]$. Once we fix elements of the finite field for the variables, say $x_1 = r_1$, $\dots$, $x_k = r_k$, the fingerprint can be evaluated in these values, and becomes a vector of matrices over the finite field. Because of the homomorphic properties presented below, the algebraic relations between polynomial fingerprints are the same as the algebraic relations between evaluated fingerprints.  

 In the next two theorems we show that the fingerprints of the results of modus ponens and substitution can be computed using the fingerprints of the inputs. The proofs are simple computations and we omit them here.
 
\begin{theorem}
Suppose that formulas $\varphi$ and $\varphi \rightarrow \psi$ have fingerprints:
$$F(\varphi) = ([\varphi], [\varphi]_{x_1}, \dots, [\varphi]_{x_k}),$$
$$F(\varphi \rightarrow \psi) = ([\varphi \rightarrow \psi], [\varphi \rightarrow \psi]_{x_1}, \dots, [\varphi \rightarrow \psi]_{x_k}).$$
Then the fingerprint of $\psi$ is:
$$F(\psi) = ([\psi], [\psi]_{x_1}, \dots, [\psi]_{x_k}),$$
where:
$$[\psi] = {A(I_2)}^{-1} \left ([\varphi \rightarrow \psi] - A(I) - A(I_1)[\varphi]\right ),  $$
$$[\psi]_{x_i} = {A(I_2)}^{-1} \left ([\varphi \rightarrow \psi]_{x_i}  - A(I_1)[\varphi]_{x_i}\right ).$$
\end{theorem} 

\begin{theorem}
    Let $\varphi$ and $\psi$ be formulas or terms. Suppose that their fingerprints are: 
    $$F(\varphi) = ([\varphi], [\varphi]_{x_1}, \dots, [\varphi]_{x_k}),$$
    $$F(\psi) = ([\psi], [\psi]_{x_1}, \dots, [\psi]_{x_k}).$$
    Let $\varphi(x_i/\psi) $ be the result of the substitution of
    $x_i$ with $\psi$ and let $X_i$ be a polynomial variable such that $A(X_i)$ is associated with the $x_i$-nodes.  Then the expression:
    $$F(\varphi(x_i/\psi)) = ([\varphi(x_i/\psi)], [\varphi(x_i/\psi)]_{x_1}, \dots, [\varphi(x_i/\psi)]_{x_k})$$
    where
    $$[\varphi(x_i/\psi)] = [\varphi] - [\varphi]_{x_i} \cdot A(X_i) + [\varphi]_{x_i} \cdot [\psi], $$ 
    and, if $j \neq i$, then:
    $$[\varphi(x_i/\psi)]_{x_j} = [\varphi]_{x_j} + [\varphi]_{x_i} [\psi]_{x_j}$$
    while if $j = i$, then:
    $$[\varphi(x_i/\psi)]_{x_i} = [\varphi]_{x_i} [\psi]_{x_i}$$
\end{theorem} 

{\bf Example}: Consider the formula
$$\varphi = (x \rightarrow y) \rightarrow (x \rightarrow z).$$
According to our definition, one has:
$$[\varphi] = [(x \rightarrow y) \rightarrow (x \rightarrow z)] = A(I) + A(I_1)[x \rightarrow y] + A(I_2) [x \rightarrow z] =$$
$$= A(I) + A(I_1)(A(I) + A(I_1)A(X) + A(I_2)A(Y)) +$$$$+ A(I_2)(A(I) + A(I_1)A(X)+A(I_2)A(Z)) = $$
$$= A(I) + A(I_1)A(I) + A(I_1)^2A(X) + $$$$+ A(I_1)A(I_2)A(Y) + A(I_2)A(I) + A(I_2)A(I_1)A(X) + A(I_2)^2 A(Z).$$
Also, one has:
$$[\varphi]_x = A(I_1)^2 + A(I_2)A(I_1),$$
$$[\varphi]_y = A(I_1)A(I_2),$$
$$[\varphi]_z = A(I_2)^2.$$
Finally, the fingerprint of $\varphi$ is:
$$F(\varphi) = ([\varphi], [\varphi]_x, [\varphi]_y, [\varphi]_z).
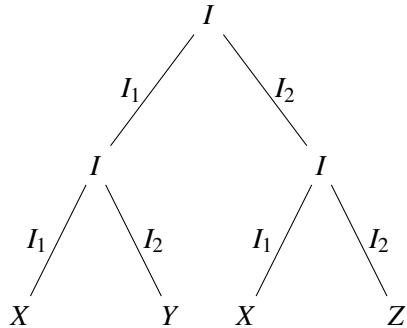
$$
\begin{center}
\begin{tikzpicture}[level distance=20mm,
    level 1/.style={sibling distance=30mm},
    level 2/.style={sibling distance=20mm},
    level 3/.style={sibling distance=10mm}]
    \node {$I$}
      child {
        node {$I$}
          child {
            node {$X$}
            edge from parent
              node[left] {$I_1$}
          }
          child {
            node {$Y$}
            edge from parent
              node[right] {$I_2$}
          }
        edge from parent
          node[left] {$I_1$}
      }
      child {
        node {$I$}
          child {
            node {$X$}
            edge from parent
              node[left] {$I_1$}
          }
          child {
            node {$Z$}
            edge from parent
              node[right] {$I_2$}
          }
        edge from parent
          node[right] {$I_2$}
      };
  \end{tikzpicture}
\end{center}


\section{Example of a formalized proof}\label{example}

In this section we apply fingerprinting for a formal mathematical proof. We consider the following framework of propositional logic. There are three axiom schemes as follows:
$$ K(\alpha, \beta):\,\,\,\,\alpha \rightarrow (\beta \rightarrow \alpha),$$
$$S(\alpha, \beta, \gamma):\,\,\,\,(\alpha \rightarrow (\beta \rightarrow \gamma)) \rightarrow ((\alpha \rightarrow \beta) \rightarrow (\alpha \rightarrow \gamma)),$$
$$N(\alpha, \beta):\,\,\,\,(\neg \alpha \rightarrow \neg \beta) \rightarrow (\beta \rightarrow \alpha).$$

We consider the following theorem:

\begin{theorem}
$$A \rightarrow A.$$
\end{theorem} 

As we will see, this particular theorem of propositional logic doesn't even need the axiom $N(\alpha, \beta)$ for its proof. However, this axiom has been mentioned for completeness. 

{\bf Classic Proof}: Consider the formula $B:= A\rightarrow A$. By making corresponding substitutions, we write:
$$S(A, B, A):\,\,\,\,(A \rightarrow (B \rightarrow A)) \rightarrow ((A \rightarrow B) \rightarrow (A \rightarrow A)),$$
$$K(A,B):\,\,\,\,A \rightarrow (B \rightarrow A),$$
$$K(A,A):\,\,\,\,A \rightarrow (A \rightarrow A),$$
At this point we observe that $K(A,A)$ is in fact:
$$K(A,A):\,\,\,\,A\rightarrow B,$$
$$MP(K(A,B), S(A,B,A)) = C:\,\,\,\, (A \rightarrow B) \rightarrow (A \rightarrow A),$$
$$MP(K(A,A), C):\,\,\,\,A \rightarrow A.$$
We also observe that the conclusion is the same as $B$. \qed 

{\bf Fingerprint translation}: We consider numeric variables $a, i, i_1, i_2$. In computations, they will be replaced by fixed choices of field elements. The corresponding matrices are:
 $$A = \begin{pmatrix}
     a & 1 \\ 0 & 1
 \end{pmatrix},\,\,\,\,
I = \begin{pmatrix}
     i & 1 \\ 0 & 1
 \end{pmatrix},\,\,\,\,
I_1 = \begin{pmatrix}
     i_1 & 1 \\ 0 & 1
 \end{pmatrix},\,\,\,\, 
I_2 = \begin{pmatrix}
     i_2 & 1 \\ 0 & 1
 \end{pmatrix},$$ 
$$ 
I_2 ^{-1} = \begin{pmatrix}
     i_2^{-1} & - i_2^{-1} \\ 0 & 1
 \end{pmatrix},$$ 
Further we define:
$$[A] = A,$$
$$[B] = I + I_1[ A ] + I_2 [A],$$
$$[B \rightarrow A] = I + I_1[B] + I_2 [A],$$
$$[K(A,B)] = I + I_1 [ A] + I_2 [B\rightarrow A] ,$$
$$[K(A,A)] = I + I_1 [A]+ I_2 [B],$$
$$[C] = I + I_1 [K(A,A)] + I_2 [B],$$
$$[S(A,B,A)] = I + I_1 [K(A,B)] + I_2 [C],$$
$$ I_2^{-1} ( [S(A,B,A)] - I - I_1[ K(A,B) ] ) = [ C ],$$
$$ I_2^{-1} ( [C] - I - I_1[ K(A,A) ] ) = [ B ].$$
\qed 

\section{First-order logic, and other logic systems}\label{firstorder}

This section presents some further constructions. Its goal is to show that the ideas presented so far - centered mainly  on modus ponens and substitution - are sufficient for a large family of logics and proof systems. 

We start with first-order logic. Here we list some general considerations for translating every first-order formula in a multivariate polynomial matrix with coefficients in $\mathbb Z$, respectively to perform fast proof verification by polynomial evaluation. An example of proof-system for first-order logic is presented in Prestel's book \cite{Prestel}. 
\begin{enumerate}
\item For building first-order terms, we follow the same tree-construction as before. Every constant $c$ from the first-order language is associated a polynomial variable $C$, respectively with the matrix:
$$
A(C) = \begin{pmatrix}
C & 1 \\ 0 & 1
\end{pmatrix}.
$$
So, if one makes the construction for the language associated with the theory of fields, to the constants $0$ and $1$ we associate polynomial variables $Z$ and $O$ and corresponding matrices $A(Z)$ and $A(O)$. The same for every first-order variable: To the variable $x$ we associate a polynomial variable $X$ and the corresponding matrix $A(X)$. If a function symbol $F$ has arity $k$, one introduces polynomial variables $F, F_1, \dots, F_k$. If the already built first-order terms $t_1, \dots, t_k$ are associated with matrices $T_1, \dots, T_k$, then the new term $F(t_1, \dots, t_k)$ is uniquely associated with the polynomial matrix $A(F) + A(F_1) T_1 + A(F_2) T_2 + \dots + A(F_k)T_k$. Suppose that we consider the additive term $x+1$ in the usual first order language of rings. This term is represented by the matrix $A(P) + A(P_1) A(X) + A(P_2)A(O)$, where $A, A_1, A_2$ are the variables associated with the operation $+$ of addition, $X$ is the polynomial variable associated with the first-order variable $x$ and $O$ is the polynomial variable associated with the constant $1$. 

\item  If a relation symbol $R$ has arity $k$, one introduces polynomial variables $R, R_1, \dots, R_k$. The atomic formula $R(t_1, \dots, t_k) $ is uniquely associated with the polynomial matrix  $A(R) + A(R_1) T_1 + A(R_2) T_2 + \dots + A(R_k)T_k$, where again the terms $t_i$ are uniquely represented by the polynomial matrices $T_i$, respectively. Also, for the equality, which is in first order logic a primitive relation, and belongs tacitely to every first-order language, we introduce unique polynomial variables $E, E_1, E_2$, such that if terms $t_1$ and $t_2$ are uniquely represented by the matrices $T_1$ and $T_2$, the atomic formula $t_1 = t_2$ is represented by the matrix $A(E) + A(E_1)T_1 + A(E_2) T_2$. 

\item The transit from (atomic) formulae to boolean combinations of (atomic) formulae is done exactly like in the propositional logic. Suppose that the first-order formulae $\varphi$ and $\psi$ are already associated with polynomial matrices $\Phi$ and $\Psi$. Then the formula $\varphi \rightarrow \psi$ is uniquely associated with the matrix $A(I) + A(I_1) \Phi + A(I_2) \Psi$, where $I, I_1, I_2$ are associated with the logic operator of implication. One deals similarly with negation, disjunction and conjunction. 

\item  We associate the quantifier $\forall\, x$ with a new polynomial variable $\bar X$, which has as unique function to represent this quantifier. If the formula $\varphi$ is uniquely associated with the matrix $\Phi$, then the formula $\forall\,x\,\,\varphi$ is uniquely associated with the matrix $A(\bar X) \Phi$. Also, the existential quantifier $\exists \, x$ is associated with a new polynomial variable $\bar {\bar X}$, such that the formula $\exists\,x\,\,\varphi$ is associated with $A(\bar {\bar X})\Phi$. We stress the fact that $X$, $\bar X$ and $\bar{\bar X}$ are three different polynomial variables, and that their role in any proof is only to represent the logical variable $x$, the quantifier $\forall\,x$, respectively the quantifier $\exists\,x$. 
\end{enumerate}

Some examples about using quantifiers. It is clear that axioms containing quantifiers, like:
$$\forall\,x\,\,\varphi \rightarrow \varphi(x/t)$$
where the term $t$ does not contain $x$, respectively
$$\forall \,x\,\,(\varphi \rightarrow \psi) \rightarrow (\phi \rightarrow \forall \,x\,\,\psi)$$
where $x$ is not a free variable of $\psi$, can be written down, and corresponding polynomial matrices can be easily computed.  Also the rule of generalization, which from $\varphi$ derives $\forall\,x\,\,\varphi$, is very easy to perform, as it consists only in multiplying the matrix $A(\bar{\bar X})$ introduced above from the left side with the corresponding matrix $\Phi$. 

We need some attention with the fact that multiplying with the matrix corresponding to a quantifier is not a commutative operation, while in logic, similar quantifiers commute inside of a block of alike quantifiers. In order to make a clean transcription of such steps, we need axioms like:
$$\forall x \forall y \forall z \,\,\varphi \rightarrow \forall y \forall z \forall x \,\,\varphi.$$
This problem is not so significant, as it appears at the first sight. For example, the commutativity of the conjunction, like the commutativity of the disjunction, needs an axiom step and a modus-ponens step, to be applied. However, it is important to use a non-commutative algebra even if we work only with propositional expressions containing only disjunctions and conjunctions, as we want expressions like $(A \wedge X) \vee (\neg X \wedge B)$ and $(A \wedge \neg X) \vee ( X \wedge B)$ to be associated with different fingerprints. 

Further difficulties in multi-sorted first order logic, modal logic  or matching  logic are treated in the same way. One has to introduce new polynomial variables for every kind of constant, variable, no matters if interpreted as element or as set, non-logical symbol, logical symbol, operator, quantified element variable or quantified set variable. In the case of operators or other symbols (like relation symbols and function symbols) of arity $k$, one always introduces $k+1$ new polynomial variables. In all these situations, every proof step can be done by a modus ponens, by a substitution, respectively by a substitution combined with the frontal multiplication with some $A(\bar X)^{-1}$, corresponding to cutting off a quantifier.
More exactly:

{\bf Multi-sorted  first-order logic.} It depends on the given syntax. But one can always replace sorts with unary predicates and modify axioms before writing down proofs, such that one finally deals with an usual first-order proof. 

{\bf Modal logic and related logics.} One can handle the modal operators box $\Box$ and diamond $\Diamond$ exactly like the quantifiers in first-order logic. One needs just one polynomial variable $B$ for the box $\Box$ and a polynomial variable $D$ for the diamond $\Diamond$. Introducing Box means a multiplication with the matrix $A(B)$ from the left side. Introducing Diamond means a multiplication with the matrix $A(D)$ from the left side. 

{\bf Matching logic.} Patterns are handled like formulas. There is no need to handle element variables differently from set variables, and the same is true for the quantified element- or set variables.  The application of patterns $\phi_1 \cdot \phi_2$ is a usual binary operation, and gets three polynomial variables. The least fix-point operator $\mu \,X. \phi$ is handled like a quantifier. This means that for a given set-variable $X$ one associates a unique polynomial variable $\bar X$ with the meaning $\mu\,X$ and the polynomial matrix $A(\bar X) \Phi$ translates  $\mu \,X. \phi$, provided that the polynomial matrix $\Phi$ translates $\phi$.

\section*{Acknowledgment}

The author had many discussions on this topic with Grigore Roșu, Xiaohong Chen, Brandon Moore, Miruna Roșca, Ruxandra Olimid-Nencioni and Traian Șerbănuță. This research was partially supported by a sponsorship from Pi Squared Inc.

\end{document}